\documentclass[12pt]{article}
\usepackage{latexsym, amssymb, amsmath}
\usepackage[all]{xy}
\usepackage[numbers,sort&compress]{natbib}
\newtheorem{theorem}{Theorem}[section]
\newtheorem{corollary}{Corollary}[section]
\newtheorem{lemma}{Lemma}[section]

\newtheorem{remark}{Remark}[section]

\def\Proof.{{\bf{Proof. }}}

\begin{document}

\title{\bf Weak convergence theorems for  a symmetric generalized hybrid mapping and an equilibrium problem}


\author{Bui Van Dinh\footnote{Faculty of Information Technology, Le Quy Don Technical University, Hanoi, Vietnam, Email: vandinhb@gmail.com}, Nguyen Ngoc Hai\footnote{
Department of Scientific Fundamentals, Vietnam Trade Union University, Hanoi, Vietnam, Email: hainn@dhcd.edu.vn}, and Do Sang Kim\footnote{Corresponding Author. Department of Applied Mathematics, Pukyong National University, Busan, Korea, Email: dskim@pknu.ac.kr}}
\maketitle

\begin{abstract}
In this paper, we introduce three new iterative methods for finding a common point of the set of fixed points of a symmetric generalized hybrid mapping and the set of solutions of an equilibrium problem in a real Hilbert space. Each method can be considered as  an combination of Ishikawa's process with the proximal point algorithm, the extragradient algorithm with or without linesearch. Under certain conditions on parameters, the iteration sequences generated by the proposed methods are proved to be weakly convergent to a solution of the problem. These results extend the previous results given in the literature. A numerical example is also provided to illustrate the proposed algorithms.
\end{abstract}
\bigskip

{\bf Keywords: }{\small  fixed point problem, equilibrium problem,  general monotonicity, \ \ \ \

 \indent extragradient method, Armijo linesearch, weak convergence.}

\bigskip



\section{Introduction}
\small
 Let $\mathbb{H}$ be a real Hilbert space endowed with an inner product $\langle \cdot , \cdot\rangle$ and the induced norm $\| \cdot \|$. We write  `$x^k \to x $', or `$x^k \rightharpoonup x $ ' iff $x^k$ converges strongly or weakly to $x$, respectively. Let $C$ be a nonempty closed convex subset of $\mathbb{H}$ and $f : C \times C \to \mathbb{R}$ be a bifunction such that
   $f(x, x) = 0$ for every $x \in C$. Such a bifunction is called an equilibrium bifunction. The equilibrium problem,  in the sense of Blum, Muu and Oettli \cite{BO,MO} (shortly EP($C, f$)), is to find $ x^* \in C $ such that
\begin{equation}\notag
  f(x^*, y) \geq 0, \  \forall y \in C.
\end{equation}
By $Sol(C, f)$, we denote the solution set of EP($C, f$).
Although problem EP($C, f$) has a simple formulation, it includes, as special cases, many important problems in applied mathematics: variational inequality problem, optimization problem, fixed point problem, saddle point problem, Nash equilibrium problem in noncooperative game, and others; see, for example, \cite{BCPP,BO,MO}, and the references quoted therein.

 Let us denote the set of fixed points of a mapping $T : C \to C$ by $Fix(T)$; that is, $Fix(T) = \{x \in C: Tx = x\}$.  Recall that a mapping $T : C \to \mathbb{H}$ is called symmetric generalized hybrid \cite{HST,KT,TWY} if there exist $\alpha, \beta, \gamma, \delta \in \mathbb{R}$ such that
$$ \alpha\|Tx - Ty\|^2 + \beta\big (\|x - Ty\|^2 + \|y - Tx \|^2 \big) + \gamma\|x - y\|^2 + \delta\big (\|x - Tx\|^2 + \|y - Ty \|^2 \big) \leq 0, \forall x, y \in C.$$ Such a mapping is called an $(\alpha, \beta, \gamma, \delta)$-symmetric generalized hybrid mapping. It is clear that
\begin{itemize}
  \item[*] (1, 0, -1, 0)-symmetric generalized hybrid mapping is nonexpansive, i.e., \\
   \centerline{ $\|Tx - Ty\| \leq \|x - y\|, \  \forall x, \  y \in C;$}
  \item[*] (2, -1, 0, 0)-symmetric generalized hybrid mapping is nonspreading, i.e., \\
  \centerline{$2\|Tx - Ty\|^2 \leq \|Tx - y\|^2 + \|Ty - x\|^2, \  \forall x, \  y \in C$ (see \cite{KoT});}
  \item[*]   (3, -1, -1, 0)-symmetric generalized hybrid mapping is hybrid, i.e., \\
  \centerline{$3\|Tx - Ty\|^2 \leq \|x - y\|^2 + \|Tx - y\|^2 + \|Ty - x\|^2, \  \forall x, \  y \in C$, (see \cite{Tak}).}
\end{itemize}
Remember that a mapping $T$ is said to be pseudo-contractive if for  all $x, \  y \in C$ and $\tau > 0$,
 $$\|x - y\| \leq \|(1+\tau)(x - y) - \tau(Tx - Ty)\|.$$
  While if $Fix(T)$ is nonempty and $\|Tx - p \| \leq \| x - p \|, \ \forall x \in C, \  p \in Fix(T)$, then $T$ is called quasi-nonexpansive. It is well-known that $Fix(T)$ is closed and convex when $T$ is quasi-nonexpansive \cite{IT}.

Most of the methods used in the literature for finding a fixed point of a mapping $T$ are derived from Mann's iteration algorithm \cite{Man}. The sequence $\{x^k\}$ is defined by the following
\begin{equation*}
\begin{cases}
x^0 \in C,\\
x^{k+1} =  \alpha_k x^k + (1 - \alpha_k)Tx^k,
\end{cases}
\end{equation*}
and under certain conditions imposed on $\{\alpha_k\}$, was proved to  converge weakly to a point in $Fix(T)$.

On the other hand,  many methods devoted to solving a monotone equilibrium problem use the proximal point algorithm:
\begin{equation}\label{1.2}
\begin{cases}
x^0 \in C,\\
\text{ find  } x^{k+1} \in C \text{ such that } f(x^{k+1}, y ) + \frac{1}{r_k}\langle y - x^{k+1}, x^{k+1} - x^k \rangle \geq 0, \forall y \in C,
\end{cases}
\end{equation}
where the sequence $\{r_k\} \subset (0, +\infty)$ and $\liminf_{k \to \infty} r_k > 0$. It was shown that the sequence $\{x^k\}$ generated by \eqref{1.2} converges weakly to a solution of $EP(C, f)$ \cite{Mou}.

Finding common elements of the solution set of an equilibrium problem and the fixed point set of a nonexpansive mapping is a task arising frequently in various areas of mathematical sciences, engineering, and economy. The motivation for studying such a problem is its possible application to mathematical models whose constraints can be expressed as fixed-point problems and/or equilbrium problems. This happens, in particular, in the practical problems as signal processing, network resource allocation, image recovery and Nash-Cournot oligopolistic equilibrium models in economy \cite{Iid,HMA}. 

For obtaining a common element of the set of fixed points of a nonexpansive mapping $T$ and the solution set of a monotone equilibrium problem EP($C, f$), Tada and Takahashi \cite{TT} proposed to combine Mann's iterative scheme with the proximal point algorithm. More precisely, the iterates $x^k, \  u^k$ are calculated as follows:
\begin{equation}\label{1.3}
\begin{cases}
x^0 \in C,\\
u^k \in C \text{ such that } f(u^k, y) + \frac{1}{r_k}\langle y - u^k, u^k - x^k \rangle \geq 0, \forall y \in C,\\
x^{k+1} =  \alpha_k x^k + (1 - \alpha_k)Tu^k.
\end{cases}
\end{equation}
The sequence $\{x^k\}$ generated by \eqref{1.3} converges weakly to some $p \in Sol(C, f) \cap Fix(T)$ provided that $\{\alpha_k\} \subset [a, b] \subset (0, 1)$ and $r_k \geq \underbar{r} > 0, \forall k$ (see, \cite[Theorem 4.1]{TT}).

Another fundamental method to find a fixed point of a mapping $T$ is Ishikawa's iteration algorithm \cite{Ish}, that is
\begin{equation}\label{1.4}
\begin{cases}
x^0 \in C,\\
y^k =  \alpha_k x^k + (1 - \alpha_k)Tx^k, \\
x^{k+1} = \beta_k x^k + (1 - \beta_k)Ty^k.
\end{cases}
\end{equation}
It was proved in \cite{Ish} that if $T$ is Lipschitzian pseudocontractive map and  $0 \leq \alpha_k \leq \beta_k \leq 1$ for all $k$, $\lim_{k \to \infty} \beta_k = 1$, $\sum_{k=1}^{\infty}(1 - \alpha_k)(1 - \beta_k) = +\infty$, then $\{x^k\}$ generated by (\ref{1.4}) converges weakly to a fixed point of mapping $T$ (see also \cite{GL}).

Motivated by these facts and recent works \cite{MA,YX,DK}, in this paper, we combine Ishikawa's algorithm with solution methods for equilibrium problems for finding a common element of the set of fixed points of a generalized hybrid mapping and the set of solutions of an equilibrium problem  in a real Hilbert space  in which the mapping $T$ is symmetric generalized hybrid, and the bifunction $f$ is monotone on $C$ or pseudomonotone on $C$ with respect to its solution set. More precisely, we propose to use the Ishikawa's algorithm for finding a fixed point of the mapping $T$ by  incorporating it with the proximal point algorithm and the extragradient algorithms with or without linesearch \cite{Kor} for solving the equilibrium problem EP($C, f$) (see also \cite{DM,DHM,FP,Kon,MQH} for more details on the extragradient algorithms). The sequences generated by the proposed algorithms are proved to converge weakly to a common solution of the symmetric generalized hybrid mapping and the equilibrium problem.

The paper is organized as follows. The next section contains some preliminaries on the metric projection, equilibrium problems and symmetric generalized hybrid mappings. The main result section is devoted to presentation of  three algorithms and their convergence in which the first one is a proximal point algorithm, the second one is  an extragradient algorithm and the last one is an extragradient algorithm with linesearch. An example and preliminary computation results are also reported.

\section{Preliminaries}
\small
Let $\mathbb{H}$ be a real Hilbert space and let $C$ be a nonempty closed convex subset of $\mathbb{H}$.
By $P_C$, we denote the metric projection operator onto $C$, that is
\begin{equation*}
 P_C(x) \in C: \Vert x - P_C(x)\Vert \leq  \Vert x - y \Vert, \  \forall y \in C.
\end{equation*}
 The following well known results on the projection operator onto a closed convex
 set will be used in the sequel.
\begin{lemma}\label{LPP} Suppose that $C$ is a nonempty closed convex subset in $\mathbb{H}$. Then
\begin{itemize}
\item[(a)] $P_C(x)$ is singleton and well defined for every $x$;

\item[(b)] $z = P_C(x)$ if and only if  $ \langle x - z  , y - z  \rangle \leq 0, \forall y \in C.$
\end{itemize}
 \end{lemma}
\begin{lemma}[Opial's condition]\label{LOP}
For any sequence $\{x^k\}\subset \mathbb{H}$ with $x^k\rightharpoonup x$, the inequality
$$\liminf\limits_{k\longrightarrow+\infty}\Vert x^k-x \Vert < \liminf\limits_{k \longrightarrow +\infty} \Vert x^k-y \Vert$$
holds for each $y\in \mathbb{H}$ with $y \neq x$.
\end{lemma}
The following result was in \cite{Zei}, page 484 (see also \cite{Sch})
\begin{lemma}\label{L1.2}
Let $\{\alpha_k\}$ be a sequence of real numbers such that $0 < a \leq \alpha_k \leq b < 1$ for all $k \in \mathbb{N}$. Let $\{v^k\}$ and $\{w^k\}$ be sequences of $\mathbb{H}$ such that, for some $c$
$$ \limsup_{k \to \infty}\|v^k\| \leq c, \;  \limsup_{k \to \infty}\|w^k\| \leq c, \text{ and } \lim_{k \to \infty}\|\alpha_kv^k + (1 - \alpha_k)w^k\| = c.$$
Then $\lim_{k \to \infty}\|v^k - w^k\| = 0.$
\end{lemma}
\begin{lemma}(\cite{TTo})\label{LTT}
Let $S$ be a nonempty closed convex subset of $\mathbb{H}$. Let $\{x^k\}$ be a sequence in $\mathbb{H}$. Suppose that, for all $p \in S$,
$$ \|x^{k+1} - p \| \leq \|x^k - p\|, \text{ for every  }  k \ = \ 0, \ 1, \ 2, \ . \ . \ .. $$
Then, $\{P_S(x^k)\}$ converges strongly to some $x^* \in S$.
\end{lemma}
In the sequel, we need the following blanket assumptions: \\

\noindent{\bf Assumptions.}
\begin{itemize}
\item[$(A_1)$] $f$ is monotone on $C$, i.e., $f(x, y) + f(y, x) \leq 0$ for all $x, \ y \in C$; 
\item[$(A_{1bis})$] $f$ is pseudomonotone on $C$ with respect to $Sol(C,f)$, i.e., $f(x, x^*)\leq 0$ for all $x \in C$, $x^* \in Sol(C, f)$;
\item[$(A_2)$]  $f(x,  \cdot )$ is convex, lower semicontinuous, and subdifferentiable on $C$, for all $x \in C$;
\item[$(A_3)$] for each $x, y, z \in C$,
$$ \limsup_{t \downarrow 0} f(tz + (1-t)x,y) \leq f(x, y); $$
 \item[$(A_{4})$] $f$ is Lipschitz-type continuous on $C$ with constants $L_1 >0$ and $L_2 >0$, i.e.,
$$f(x, y)+f(y, z) \geq f(x, z) - L_1 \Vert x-y \Vert^2 - L_2 \Vert y-z \Vert^2, \; \forall x, y, z\in C;$$
\item[$(A_{4bis})$] $f$ is jointly weakly continuous on $C\times C$ in the sense that if $x, y \in C$ and $\{x^k\}$,  $ \{y^k\} $ are two sequences in $C$
\indent converging weakly to $x$ and $y$ respectively, then $f(x^k, y^k)$ converges to $f(x, y)$;
\item[$(A_5)$] $T$ is a ($\alpha, \beta, \gamma, \delta$)-symmetric generalized hybrid self-mapping of $C$ such that $(1) \ \alpha + 2\beta + \gamma \geq 0$,
\indent $(2) \ \alpha + \beta > 0$,  $(3) \ \delta \geq 0 $, and $Fix(T)$ is nonempty.
 \end{itemize}
The following lemma is well-known in theory of monotone equilibrium problems.
\begin{lemma}(\cite{BO,CH})\label{CH}
For $\rho > 0$, $x \in \mathbb{H}$, define a mapping  $T_{\rho}: \mathbb{H} \to C$ as follows
\begin{equation*}\label{}
T_\rho(x)=\Big\{u \in C: f(u, y)+\displaystyle\frac{1}{\rho}\langle y - u, u - x \rangle \geq 0 \; \; \forall y \in C \Big\}.
\end{equation*}
Then under assumptions ($A_1$), ($A_2$), and ($A_3$) the following statements hold: \\
\indent (i) $T_\rho$ is well defined and single-valued;\\
\indent (ii) $T_\rho$ is firmly nonexpansive, i.e., for any $x, y \in \mathbb{H}$,
$$\Vert T_\rho (x) -T_\rho (y) \Vert^2 \leq \langle T_\rho(x) - T_\rho (y), x - y \rangle;$$
\indent (iii) Fix$(T_\rho) = \text{Sol}(C, f)$;\\
\indent (iv) Sol$(C, f)$ is closed and convex.
\end{lemma}
For each   $z, \  x \in C$, by $\partial_2 f(z , x)$ we denote the subdifferential of the convex
  function $f(z, .)$ at $x$, i.e.,
$$ \partial_2 f(z , x) := \{ w \in \mathbb{H} : f(z, y) \geq f(z, x) + \langle w, y - x \rangle, \ \forall y \in C \}. $$
In particular,
$$ \partial_2 f(z , z)  = \{ w \in \mathbb{H} : f(z, y) \geq  \langle w, y - z \rangle, \ \forall y  \in C \}.$$
Let $\Omega$ be an open convex set containing $C$. The next lemma can be considered as an infinite-dimensional version of Theorem 24.5 in \cite{Roc}
\begin{lemma}\label{Lem2.2}\cite[Proposition 4.3]{VSH}
Let $f : \Omega \times \Omega \to \mathbb{R}$ be a function satisfying conditions $(A_2)$ on $C$ and ($A_{4bis}$) on $\Omega$. Let $\bar{x}, \bar{y} \in \Omega$ and $\{x^k\}$, $\{y^k\}$ be two sequences in $\Omega$ converging weakly to $\bar{x}, \bar{y}$, respectively. Then, for any $\epsilon > 0$, there exist $\eta >0$ and $k_{\epsilon} \in \mathbb{N}$ such that
$$ \partial_2 f(x^k, y^k) \subset \partial_2 f(\bar{x}, \bar{y}) + \frac{\epsilon}{\eta}B, $$
for every $k \geq k_\epsilon$, where $B$ denotes the closed unit ball in $\mathbb{H}.$
\end{lemma}
\begin{lemma}\label{Lem2.6}  Let the bifunction $f$ satisfy the assumptions $(A_2)$ on $C$ and (${A}_{4bis}$) on $\Omega$, and $\{x^k \} \subset C $, $0 < \underbar{$\rho$} \leq \bar{\rho} $,  $ \{\rho_k\} \subset [\underbar{$\rho$} , \ \bar{\rho}] $. Consider the sequence $\{y^k\}$ defined as follows
$$ y^k = \arg\min\Big\{  f(x^k, y) +  \frac{1}{2\rho_k} \| y-x^k\|^2: \ y\in C \Big\}.$$
If $\{x^k \}$ is bounded, then $\{y^k\}$ is also bounded.
  \end{lemma}
{\bf Proof.}
Firstly, we show that if $\{x^k \}$  converges weakly to $x^*$, then $\{y^k\}$ is bounded. \\
Since
$$ y^k = \text{arg}\min\Big\{  f(x^k, y) +  \frac{1}{2\rho_k} \| y-x^k\|^2: \ y\in C \Big\}, \forall k$$
and since
$$ f(x^k, x^k) +  \frac{1}{2\rho_k} \| x^k-x^k\|^2 = 0, \forall k,$$
we obtain
$$ f(x^k, y^k) +  \frac{1}{2\rho_k} \| y^k-x^k\|^2 \leq 0, \ \forall k. $$
In addition, for all $w^k \in \partial_2f(x^k, x^k)$ we have
 $$ f(x^k, y^k) +  \frac{1}{2\rho_k} \| y^k-x^k\|^2  \geq \langle w^k, y^k - x^k \rangle + \frac{1}{2\rho_k} \|y^k - x^k\|^2.$$
This implies $$ -\|w^k\| \|y^k - x^k \| + \frac{1}{2\rho_k} \|y^k - x^k\|^2 \leq 0.$$ Hence,
$$\|y^k - x^k\| \leq 2\rho_k \|w^k\|, \ \forall k.$$
 Since $\{\rho_k\}$ is bounded,   $\{x^k\}$ converges weakly to $x^*$ and $w^k \in \partial_2f(x^k, x^k)$, it follows from Lemma~\ref{Lem2.2} that the sequence $\{w^k\}$ is bounded. The sequence $\{x^k\}$ being bounded, we get that $\{y^k\}$ is also bounded.

Now we prove Lemma~\ref{Lem2.6}. Suppose that $\{y^k\}$ is unbounded, i.e., there exists a subsequence $\{y^{k_i}\} \subseteq  \{y^k\}$ such that $\lim\limits_{i \to \infty}\|y^{k_i}\| = + \infty$. By the boundedness of  $\{x^k\}$, the subsequence  $\{x^{k_i}\}$ is also bounded, and without loss of generality, we may assume that $\{x^{k_i}\}$ converges weakly to some $ x^*$. By the same argument as above, we obtain that $\{y^{k_i}\}$ is bounded, which contradicts to the fact that $\lim\limits_{i \to \infty}\|y^{k_i}\| = + \infty$. 
Therefore  $\{y^k\}$ is bounded. \hfill$\Box$

The following lemmas give us a characterization of fixed point set of symmetric generalized hybrid mappings.

\begin{lemma}\label{L2.5}\cite{KT}
Let $C$ be a nonempty closed convex subset of  $\mathbb{H}$ . Assume that $T$ is an
($\alpha, \beta, \gamma, \delta$)-symmetric generalized hybrid self-mapping of $C$ such that $Fix(T) \neq \emptyset$ and the conditions $(1) \ \alpha + 2\beta + \gamma \geq 0$,  $(2) \ \alpha + \beta > 0$ and $(3) \ \delta \geq 0 $ hold. Then T is quasi-nonexpansive.
\end{lemma}
\begin{lemma}\label{L2.7}\cite{HST}
Let $C$ be a nonempty closed convex subset of  $\mathbb{H}$ . Assume that $T$ is an
($\alpha, \beta, \gamma, \delta$)-symmetric generalized hybrid self-mapping of $C$ such that $Fix(T) \neq \emptyset$ and the conditions $(1) \ \alpha + 2\beta + \gamma \geq 0$,  $(2) \ \alpha + \beta > 0$ and $(3) \  \delta \geq 0 $ hold. Then $I - T$ is demiclosed at $0$, i.e., $x^k \rightharpoonup \bar{x}$ and $x^k - Tx^k \to 0$ imply $\bar{x} \in Fix(T)$.
\end{lemma}

\section{Main Results}
\small
Based on the idea of Tada and Takahashi \cite{TT}, we now combine algorithms for equilibrium problems and Ishikawa's process for fixed point problems to get the following algorithms.  \\

 \noindent{\bf Algorithm 1}
\begin{itemize}
   \item[]{\bf Initialization.} Pick  $x^0 \in C$, choose parameters  $  \b{$\beta$} \in (0, 1)$, $  \b{$\rho$} > 0$;  $\{ \rho_k \} \subset [ \b{$\rho$}, \  + \infty)$, $ \{ \alpha_k \} \subset [0, 1], \ \lim_{k \to \infty} \alpha_k = 1$;  $\{\beta_k\} \subset [\b{$\beta$}, \bar{\beta}] \subset [0, 1)$.
   \item[]{\bf Iteration $k$} (k = 0, 1, 2, ...). Having $x^k$  do the following steps:
	\begin{itemize}
	 \item[]{\it Step 1.} Find $u^k \in C$ such that
            $$f(u^k, y) +  \frac{1}{\rho_k} \langle y - u^k, u^k - x^k \rangle \geq 0, \text{ for all } y \in C$$
   	 \item[]{\it Step 2.} Compute
                   \begin{align*}
                      v^{k} &= \alpha_k x^k + (1 -\alpha_k) Tx^k,\\
                     x^{k+1} &= \beta_k v^k + (1 -\beta_k) Tu^k,
                      \end{align*}
	 and go to Step 1  with $k$ is replaced by  $k +1 $.\\
\end{itemize}
\end{itemize}

From Lemma~\ref{CH}, we have that $\{u^k\}$ is well defined. Hence $\{x^k\}$ is well defined. The following theorem establishes the convergence of Algorithm 1.

\begin{theorem}\label{The3.1} Suppose that the set $S = Sol(C, f) \cap Fix(T)$ is nonempty. Then  under assumptions ($A_1$), ($A_2$), $(A_3)$, and $(A_5)$, the sequences $\{x^k\}$, $\{u^k\}$ generated by Algorithm 1 converge weakly to $x^* \in S$, where $x^* = \lim_{k \to \infty} P_S (x^k)$.
      \end{theorem}
{\bf Proof.} Let $k$ be fixed. By  definition of  $u^k$, we can write $u^k = T_{\rho_k}(x^k)$. Taking some  $q \in \ S$, i.e., $q \in  Sol(C, f) \cap Fix(T)$ and using the non-expansiveness of $T_{\rho_k}$,  we have
\begin{equation}\label{3.3}
   \|u^k - q \| = \| T_{\rho_k}(x^k) - T_{\rho_k}(q) \leq \|x^k - q \|.
\end{equation}
From Step 2, we have
\begin{equation*}
\begin{aligned}
\|v^k - q\| & = \| \alpha_kx^k + (1 - \alpha_k) Tx^k - q \| \\
                & \leq \alpha_k\|x^k - q\| + (1 - \alpha_k)\|Tx^k - q \|.
\end{aligned}
\end{equation*}
Since $T$ is a ($\alpha, \beta, \gamma, \delta)$-symmetric generalized hybrid mapping with $\alpha + 2\beta + \gamma \geq 0$, $\alpha + \beta > 0$, $\delta \geq 0$, it follows from Lemma~\ref{L2.5} that $T$ is quasi-nonexpansive. So
\begin{equation}\label{3.4}
\| v^k - q \|  \leq  \| x^k - q \|.
\end{equation}
Similarly,
\begin{equation}\notag
\begin{aligned}
\|x^{k+1} - q\| & = \|\beta_k v^k + (1 - \beta_k) Tu^k - q \|\\
&\leq \beta_k \|v^k - q \| + (1 - \beta_k) \|Tu^k - q \|\\
&\leq \beta_k \|v^k - q \| + (1 - \beta_k) \|u^k - q \|.
\end{aligned}
\end{equation}
Combining with (\ref{3.3}) and (\ref{3.4}) yields
\begin{equation}\label{3.5}
\| x^{k+1} - q \|  \leq  \| x^k - q \|.
\end{equation}
Since (\ref{3.5}) holds for all $k$, we have that $\lim_{k \to   \infty}\| x^k -  q \| $ does exist. Let $\tau = \lim_{k \to   \infty}\| x^k -  q \|$. Consequently, the sequence $\{x^k\}$ is bounded, and from (\ref{3.3}), (\ref{3.4}), we get that $\{u^k\}$, $\{v^k\}$ are also bounded.\\
In addition, by Lemma~\ref{CH} ($ii$), we have
\begin{equation}\notag
\begin{aligned}
\|u^{k} - q \|^2 & = \|T_{\rho_k}(x^k) - T_{\rho_k}(q)\|^2  \leq \langle T_{\rho_k}(x^k) - T_{\rho_k}(q) , x^{k} - q \rangle \\
                               & = \langle u^k - q, x^k - q \rangle  = \frac{1}{2}\big[\|u^k - q\|^2 + \|x^{k} - q \|^2 - \|u^k - x^k\|^2 \big].
\end{aligned}
\end{equation}
So
\begin{equation}\label{3.6}
\|u^k - q\|^2 \leq \|x^k - q \|^2 - \|u^k - x^k \|^2, \ \forall k, \ \forall q \in S.
\end{equation}
By definition of $x^{k+1}$, we get
\begin{equation}\notag
\begin{aligned}
\|x^{k+1} - q \|^2 & = \|\beta_k(v^{k} - q) + (1 - \beta_k)(Tu^k - q) \|^2 \\
                                       & = \beta_k\|v^{k} - q \|^2 + (1 - \beta_k)\|Tu^k - q \|^2  - \beta_k(1 - \beta_k) \|Tu^k - v^k \|^2 \\
                                       & \leq \beta_k\|v^{k} - q \|^2 + (1 - \beta_k)\|Tu^k - q \|^2  \\
                                      & \leq \beta_k\|v^{k} - q \|^2 + (1 - \beta_k)\|u^k - q \|^2.     \\
\end{aligned}
\end{equation}
In view of (\ref{3.4}) and (\ref{3.6}) we deduce that
\begin{equation}\notag
\|x^{k+1} - q \|^2 \leq \|x^{k} - q \|^2 - (1-\beta_k)\|u^{k} - x^k \|^2.
\end{equation}
Hence
\begin{equation}\label{3.8}
 (1-\bar{\beta})\|u^{k} - x^k \|^2  \leq \|x^{k} - q \|^2  - \|x^{k+1} - q \|^2.
\end{equation}
Because $\lim_{k \to \infty}\|x^k - q\| = \tau$, it follows from (\ref{3.8}) that
\begin{equation}\label{3.9}
\lim_{k \to \infty}\|u^{k} - x^k \| = 0.
\end{equation}
Since $\lim_{k \to \infty} \alpha_k = 1$, we have
\begin{equation}\label{3.10}
\lim_{k \to \infty} \|v^k - x^k \| = \lim_{k \to \infty} (1 - \alpha_k)\|x^k - T x^k \| = 0.
\end{equation}
It is clear that
\begin{equation}\notag
\lim_{k \to \infty} \|\beta_k(v^k - q) + (1-\beta_k)(Tu^k -q) \| = \lim_{k \to \infty} \|x^{k+1} - q \| = \tau.
\end{equation}
Hence
\begin{equation*}
\lim_{k \to \infty}\sup \|v^k - q \| \leq \lim_{k \to \infty}\sup\|x^k - q \| = \tau,
\end{equation*}
and
\begin{equation*}
\lim_{k \to \infty}\sup \|Tu^k - q \| \leq \lim_{k \to \infty}\sup\|u^k - q \| \leq \lim_{k \to \infty}\sup\|x^k - q \|  = \tau,
\end{equation*}
in the light of Lemma~\ref{L1.2}, yield
\begin{equation}\label{3.11}
\lim_{k \to \infty} \|v^k - Tu^k \| = 0.
\end{equation}
On the other hand, for each $k$, we can write
\begin{equation*}
\begin{aligned}
 \|Tu^k - u^k \| & \leq \| Tu^k - v^k \| + \|v^k - x^k \| +  \|x^k - u^k\|,
\end{aligned}
\end{equation*}
and combining this inequality with (\ref{3.9}), (\ref{3.10}), and  (\ref{3.11}) we can deduce in the limit that
\begin{equation}\label{3.15}
\lim_{k\to \infty} \|Tu^k - u^k \| = 0.
\end{equation}
Next we show that any weak accumulation point of $\{x^k\}$ belongs to $S$. Indeed, suppose that $\{x^{k_i}\} \subset \{x^k\}$ and $x^{k_i} \rightharpoonup x^*$  as $i \to \infty$. From (\ref{3.9}) one has $u^{k_i} \rightharpoonup x^*$ as $i \to \infty$, and
since $\lim\inf_{k \to \infty} \rho_k > 0$ that
\begin{equation}\label{3.16}
\lim_{k \to \infty}\frac{\|u^{k} - x^k \|}{\rho_k} = 0.
\end{equation}
By definition of $u^k$, we get
\begin{equation*}
 f(u^{k}, y) + \frac{1}{\rho_k} \langle y - u^k , u^k - x^k \rangle \geq 0, \ \text{ for all } y \in C,
\end{equation*}
and by monotonicity of $f$, we can write
\begin{equation*}
 \frac{1}{\rho_k} \langle y - u^k , u^k - x^k \rangle \geq f(y, u^{k}) , \ \text{ for all } y \in C.
\end{equation*}
So
\begin{equation}\label{3.17}
  \langle y - u^{k_i} , \frac{u^{k_i} - x^{k_i}}{\rho_{k_i}} \rangle \geq f(y, u^{k_i}) , \ \text{ for all } y \in C.
\end{equation}
Letting $i \to \infty$, by the continuity of $f$ and  (\ref{3.16}), we obtain in the limit from (\ref{3.17}) that
$$ 0 \geq  f(y, x^*), \text{ for all } y \in C.$$
Suppose that $t \in (0, 1]$, $y \in C$, let $y_t = ty + (1-t)x^*$. Since $y \in C$ and $x^* \in C$, it follows that $y_t \in C$ and hence $f(y_t, x^*) \leq 0.$ So, we have
$$ 0 = f(y_t, y_t) \leq tf(y_t, y) + (1-t)f(y_t, x^*) \leq tf(y_t, y).$$
Therefore
$$ f(y_t, y) \geq 0 \text{ for all } t \in (0; 1] \text{ and all } y \in C. $$
By taking the limit as $t \downarrow 0$ and using ($A_3$) we get $f(x^*, y) \geq 0, \text{ for all } y \in C,$
which  means that $x^*$ is a solution of EP($C, f$).\\
By virtue of  (\ref{3.15}), we obtain $\lim_{i \to \infty}\|Tu^{k_i} - u^{k_i}\| = 0$. Since $u^{k_i} \rightharpoonup x^*$ and $I-T$ is demiclosed at zero, by Lemma~\ref{L2.7}, we get $Tx^* = x^*$, i.e., $x^* \in Fix(T).$ \\
Therefore $x^* \in S$.

To complete the proof, we must show that the whole sequence  $\{x^k\}$ converges weakly to $x^*$. Indeed, if there exists a subsequence $\{x^{l_i}\}$ of $\{x^k\}$ such that $x^{l_i} \rightharpoonup \hat{x}$ with $ \hat{x} \neq x^* $, then we have that $ \hat{x} \in S$, and by Opial's condition, that
\begin{equation}\notag
\begin{aligned}
\liminf\limits_{i \to +\infty}\Vert x^{l_i} - \hat{x} \Vert & < \liminf\limits_{i \to  +\infty}\Vert x^{l_i} - x^*\Vert\\
&=\liminf\limits_{j \to +\infty}\Vert x^{k} - x^* \Vert\\
&=  \liminf\limits_{j \to +\infty}\Vert x^{k_j}-x^* \Vert \\
&<\liminf\limits_{j \to +\infty}\Vert x^{k_j}-\hat{x} \Vert\\
&=\liminf\limits_{i \to +\infty}\Vert x^{l_i}-\hat{x}\Vert.
\end{aligned}
\end{equation}
This is a contradiction. Hence $\{x^k\}$ converges weakly to $x^*$ and from \eqref{3.9} and \eqref{3.10}, we deduce immediately that $\{u^k\} $, $\{v^k\}$ also converge weakly to $x^*$.
From (\ref{3.5}) and Lemma~\ref{LTT} we have that $\{ P_S(x^k)\}$ strongly converges to some $\hat{x} \in S.$
In addition, from Lemma~\ref{LPP}, we derive
$$ \langle x^* - P_S(x^k), x^k - P_S(x^k) \rangle \leq 0, \forall k$$
and by taking  the limit as $k \to \infty$, the above inequality becomes
$$ \langle x^* - \hat{x} , x^* - \hat{x} \rangle \ = \ \|x^* - \hat{x}\|^2  \leq 0. $$
Therefore $x^* = \hat{x}$. This completes the proof. \hfill$\Box$ \\
When $\alpha_k = 1, \forall k$ we have that $v^k = x^k, \forall k$ and Theorem 1 becomes
\begin{corollary}\label{C1}
Suppose that the set $S = Sol(C, f) \cap Fix(T)$ is nonempty and assumptions ($A_1$), ($A_2$), $(A_3)$, and $(A_4)$ are satisfied. Consider the sequences $\{x^k\}$, $\{u^k\}$ generated by $x^0 \in \mathbb{H}$, $0 < \b{$\rho$} $;
$ \{ \rho_k \} \subset [ \b{$ \rho $}, \  + \infty) $,   $\{\beta_k\} \subset [\b{$ \beta $}, \bar{\beta}] \subset (0, 1)$ and
\begin{equation}\notag
\begin{cases}
u^k \in C \text{ such that } f(u^k, y) +  \frac{1}{\rho_k} \langle y - u^k, u^k - x^k \rangle \geq 0, \text{ for all } y \in C \\
x^{k+1} = \beta_k v^k + (1 -\beta_k) Tu^k.
\end{cases}
\end{equation}
Then, $\{x^k\}$, $\{u^k\}$ converge weakly to $x^* \in S$, where $x^* = \lim_{k \to \infty} P_S (x^k).$
\end{corollary}
\begin{remark}
Theorem 4.1 of Tada and Takahashi \cite{TT} is a special case of corollary~\ref{C1}, because nonexpansive mappings are symmetric generalized hybrid mappings.
\end{remark}
\begin{remark}
For each $x^k \in C$, $l_k(x,y) = \frac{1}{\rho_k}\langle y - x, x - x^k \rangle$ is strongly monotone on $C$ with constant $\tau = \frac{1}{\rho_k}$ (i.e., $l(x,y) + l(y,x) \leq -\tau \|x - y\|^2, \ \forall x, y \in C$). Hence, if $f$ is monotone on $C$, then the function $f_k(x,y) = f(x,y) + l_k(x,y)$ is strongly monotone with constant $\tau$, and therefore,  {\it Algorithm 1} is well defined and to find $u^k$ at {\it Step 1}, we can apply some existing methods, see, for instance \cite{BCPP,MQ}. However, if $f$ is pseudomonotone on $C$, the bifunction $f_k$ may not be strongly monotone, even not be pseudomonotone on $C$; see,  counterexample 2.1 in \cite{TYY}, example 2.8 in \cite{DM2}, so  we can not apply the available methods using the monotonicity of the bifunction $f_k$ to find $u^k$ directly.
\end{remark}

To find a solution of  pseudomonotone equilibrium problem, Tran et al. \cite{MQH} proposed to use the extragradient algorithm introduced by Korpelevich \cite{Kor} for finding saddle points and other related problems. Now, we combine the extragradient algorithm with Ishikawa process to get the following algorithm for symmetric generalized hybrid mapping and equilibrium problem.\\

 \noindent{\bf Algorithm 2.}
\begin{itemize}
   \item[]{\bf Initialization.} Pick  $x^0  \in C$, choose parameters $\{ \rho_k \} \subset [ \b{$\rho$}, \  \bar{\rho} ]$, with $0 < \b{$\rho$} \leq \bar{\rho} < \min\{\frac{1}{2L_1}, \frac{1}{2L_2}\} $, \\
        \indent $ \{ \alpha_k \} \subset [0, 1], \ \lim_{k \to \infty} \alpha_k = 1$, $\{\beta_k\} \subset [\b{$\beta$}, \bar{\beta}] \subset (0, 1)$.
   \item[]{\bf Iteration $k$} (k = 0, 1, 2, ...). Having $x^k$  do the following steps:
	\begin{itemize}
	 \item[]{\it Step 1.} Solve successively the strongly convex programs
            $$\min\Big\{  \rho_k f(x^k, y) +  \frac{1}{2} \| y-x^k\|^2: \ y\in C\Big\} \ \ \ \ \ \ \ \ \ \ \  \eqno CP(x^k)$$
            $$ \min\Big\{  \rho_k f(y^k, y) +  \frac{1}{2} \| y-x^k\|^2: \ y\in C\Big\} \ \ \ \ \ \  \eqno CP(y^k, x^k)$$
             to obtain their unique solutions $y^k$ and $z^k$ respectively.
	 \item[]{\it Step 2.} Compute
                   \begin{align*}
                      v^{k} &= \alpha_k x^k + (1 -\alpha_k) Tx^k,\\ 
                     x^{k+1} &= \beta_k v^k + (1 -\beta_k) Tz^k, 
                      \end{align*}
	 and go to Step 1  with $k$ is replaced by  $k +1 $.\\
\end{itemize}
\end{itemize}

Before proving the convergence of  this algorithm, let us recall the following result which was proved in \cite{Anh}
\begin{lemma}\label{Lem3.1}\cite{Anh}
Suppose that $f$ satisfies assumption ($A_2$) and $x^* \in \text{ Sol}(C, f)$,  then we have:
\begin{itemize}
\item[(i)]  $\rho_k[f(x^k, y)-f(x^k, y^k)] \geq \langle y^k-x^k, y^k-y \rangle,\;\forall y \in C.$
\item[(ii)]  If, in addition, $f$ satisfies assumptions $({A}_{1bis})$, and $(A_{4})$, then we have:
$$\| z^k-x^* \|^2  \leq \|x^k-x^* \|^2- (1 - 2\rho_k L_1) \| x^k-y^k \|^2 - (1-2 \rho_k L_2) \| y^k - z^k \|^2,\;\;\forall k.$$
\end{itemize}
\end{lemma}

\begin{theorem}\label{T4.1} Suppose that the bifunction $f$ and the mapping $T$ satisfy the assumptions (${A}_{1bis}$), ($A_2$), $({A}_3)$, $(A_{4})$  and $(A_5)$, respectively, and the set $S = Sol(C, f) \cap Fix(T)$ is nonempty. Then the sequences $\{x^k\}$, $\{y^k\}$, $\{z^k\}$ generated by Algorithm 2 converge weakly to $x^* \in S$, where $x^* = \lim_{k \to \infty} P_S (x^k)$.
      \end{theorem}
{\bf Proof.}
Take  any $q \in \ S$, from Lemma~\ref{Lem3.1} we have
$$\| z^k - q \|^2  \leq \|x^k - q \|^2- (1 - 2\rho_k L_1) \| x^k-y^k \|^2 - (1-2 \rho_k L_2) \| y^k - z^k \|^2,\;\;\forall k.$$
Because $ 0 < \b{$\rho$} \leq \rho_k \leq{\bar{\rho}} < \text{min}\{\frac{1}{2L_1}, \frac{1}{2L_2}\}$, we get
\begin{equation}\label{e3.3}
   \|z^k - q \| \leq \|x^k - q \|.
\end{equation}
Arguing similarly as in the proof of Theorem 1, we conclude that
\begin{equation}\label{e3.4}
\| v^k - q \|  \leq  \| x^k - q \|,
\end{equation}
and
\begin{equation}\label{e3.5}
\| x^{k+1} - q \|  \leq  \| x^k - q \|.
\end{equation}
Hence
\begin{equation}\label{e3.9}
\lim_{k \to \infty} \| x^{k} - q \| = \tau.
\end{equation}
  In view of (\ref{e3.3}) and (\ref{e3.4}), we get $\{z^k\}$, $\{v^k\}$ are also bounded.\\
We have
\begin{equation}\notag
\begin{aligned}
\|x^{k+1} - q \|^2 & = \|\beta_k(v^{k} - q) + (1 - \beta_k)(Tz^k - q \|^2 \\
                                      & \leq \beta_k\|v^{k} - q \|^2 + (1 - \beta_k)\|z^k - q \|^2.     \\
\end{aligned}
\end{equation}
Combining with (\ref{e3.4}) and Lemma~\ref{Lem3.1}, yields
\begin{equation*}
\| x^{k+1} - q \|^2 \leq   \|x^k - q\|^2 - (1 - \beta_k) \big[ (1 - 2\rho_k L_1) \| x^k-y^k \|^2 - (1-2 \rho_k L_2) \| y^k - z^k \|^2  \big].
\end{equation*}
Therefore
\begin{equation}\label{e3.11}
(1 - \beta_k)\big[(1 - 2\rho_k L_1) \| x^k-y^k \|^2 + (1-2 \rho_k L_2) \| y^k - z^k \|^2 \big]  \leq \big(\|x^k - q \| - \| x^{k+1} - q \| \big) \big( \|x^k - q \| + \| x^{k+1} - q \| \big).
\end{equation}
Since $0 < 1 -\bar{\beta} \leq 1- \beta_k$; $0 < \b{$\rho$} \leq \rho_k \leq \bar{\rho} < \min\{\frac{1}{2L_1}, \frac{1}{2L_2}\}$, and (\ref{e3.9}), we can conclude from (\ref{e3.11}) that
\begin{equation}\label{e3.12}
\lim_{k \to \infty} \|x^k - y^k \| = 0.
\end{equation}
\begin{equation}\label{e3.13}
\lim_{k \to \infty} \|y^k - z^k \| = 0.
\end{equation}
By the triangle inequality, we deduce from \eqref{e3.12} and \eqref{e3.13} that
\begin{equation}\label{e3.14}
\lim_{k \to \infty} \|x^k - z^k \| = 0.
\end{equation}
Using the same argument as in Theorem 1, we have
\begin{equation}\label{e3.18}
\lim_{k\to \infty} \|Tz^k - z^k \| = 0.
\end{equation}
Now, suppose that $\{x^{k_i}\}$ is any subsequence of $ \{x^k\}$ such that $\{x^{k_i} \}$ converges weakly to $x^*$  as $i \to \infty$. In view of (\ref{e3.12}) and (\ref{e3.14}), we obtain $y^{k_i} \rightharpoonup x^*$, and $z^{k_i} \rightharpoonup x^*$ as $i \to \infty$.\\
Replacing $k$ by $k_i$ in assertion (i) of Lemma~\ref{Lem3.1}, it yields
\begin{equation*}
\rho_{k_i}\big[ f(x^{k_i}, y) - f(x^{k_i}, y^{k_i}) \big] \geq  \langle x^{k_i} - y^{k_i} , y - y^{k_i} \rangle, \ \forall y \in C.
\end{equation*}
Hence
\begin{equation}\label{e3.19}
\rho_{k_i}\big[ f(x^{k_i}, y) - f(x^{k_i}, y^{k_i}) \big] \geq   -\|x^{k_i} - y^{k_i}\| \| y - y^{k_i}\|.
\end{equation}
Letting $i \to \infty$, by the continuity of $f$ and  (\ref{e3.12}), we obtain in the limit from (\ref{e3.19}) that
$$ f(x^*, y) - f(x^*, x^*) \geq 0. $$
So, $ f(x^*,y)  \geq 0, \ \forall y \in C, $ which  means that $x^*$ is a solution of EP($C, f$).\\
From (\ref{e3.18}), one has $\lim_{i \to \infty}\|Tz^{k_i} - z^{k_i}\| = 0$. Because $z^{k_i} \rightharpoonup x^*$ and $I - T$ is demiclose at zero, using Lemma~\ref{L2.7}, we obtain $Tx^* = x^*$, i.e., $x^* \in Fix(T).$ \\
Hence $x^* \in S$. The rest of the proof can be done similarly to Theorem 1 so we obmit it.
The proof is completed.
\hfill$\Box$

\begin{remark}
The parameters $\{\rho_k\}$  in Algorithm 2 are determined by the Lipschitz constants $L_1$ and $L_2$ of $f$. However, in general,  these constants are usually difficult to estimate or $f$  doesn't satify the Lipschitz condition, so we can not apply Algorithm 2 to solve the above problem directly.
\end{remark}
To solve equilibrium problem EP($C, f$) when $f$ doesn't satisfy Lipschitzian type conditions Tran et al \cite{MQH}, Dinh and Muu \cite{DM} introduced linesearch methods. The following algorithm can be seen as a combination of linesearch algorithm and Ishikawa's process for finding a common point of solution set of equilibrium problem and the set of fixed points of symmetric generalized hybrid mapping.\\ 

\noindent{\bf Algorithm 3}
\begin{itemize}
   \item[]{\bf Initialization.} Pick  $x^0  \in C$, choose parameters $\eta, \mu \in (0, 1); \ 0 < \b{$\rho$} \leq \bar{\rho}$, $\{ \rho_k \} \subset [ \b{$\rho$}, \ \bar{\rho}]$; $ \{ \alpha_k \} \subset [0, 1],$  $\lim_{k \to \infty} \alpha_k = 1$;    $\{\beta_k\} \subset [\b{$\beta$}, \bar{\beta}] \subset (0, 1)$; $\gamma_k \in [\b{$\gamma$}, \bar{\gamma}] \subset (0, 2)$.
\item[]{\bf Iteration $k$} (k = 0, 1, 2, ...). Having $x^k$  do the following steps:
	\begin{itemize}
	\item[]{\it Step 1.}
		$$y^k = \arg\min\Big\{  \rho_k f(x^k, y) +  \frac{1}{2} \| y-x^k\|^2: \ y\in C\Big\}  \ \ \  \ \ \ \ \ \ \ \ \      	\ \eqno $$
	\item[]{\it Step 2.} (Armijo linesearch rule) Find $m_k$ as the smallest positive integer 	   number $m$ such that
	\begin{equation}\label{4.1}
	\begin{cases}
	z^{k,m} = (1 - \eta^m)x^k + \eta^my^m \\
	f(z^{k,m}, x^k) - f(z^{k,m}, y^k)  \geq \frac{\mu}{2 \rho_k}\|x^k - y^k\|^2.
	\end{cases}
	\end{equation}
	Set $\eta_k = \eta^{m_k}$, $z^k = z^{k, m_k}$.
	\item[]{\it Step 3.} Select $w^k \in \partial_2f(z^k, x^k)$, and compute $u^k = P_C(x^k - \gamma_k \sigma_k w^k)$, \\
     where $\sigma_k = \frac{f(z^k, x^k)}{\|w^k\|^2}$.
     \item[]{\it Step 4.} Compute
      \center{	$v^{k} = \alpha_k x^k + (1 -\alpha_k) Tx^k$,} \\
 $x^{k+1} = \beta_k v^k + (1 -\beta_k) Tu^k.$ 
          \end{itemize}
\end{itemize}
To prove the convergence of Algorithm 3 we need the following lemma.
\begin{lemma}\label{Lem4.1}\cite{MQH}
Suppose that $p \in \text{ Sol}(C, f)$, then under assumptions $({A}_{1bis})$ and $(A_2)$. Then, we have:
\begin{itemize}
\item[(a)]  The linesearch is well defined;
\item[(b)]  $ f(z^k, x^k) > 0$;
\item[(c)] $0 \not\in \partial_2f(z^k, x^k)$;
\item[(d)] \begin{equation}\notag
                   \|u^k - p\| \leq \|x^k - p\|^2 - \gamma_k( 2 - \gamma_k)(\sigma_k\|w^k\|)^2.
                 \end{equation}
\end{itemize}
\end{lemma}
{\bf Proof.} The proof of  Lemma~\ref{Lem4.1} when $\mathbb{H}$ is a finite dimensional space could be found, for instance \cite{MQH}. When its dimension is infinite, it can be done by the same way. So we omit it.
\begin{theorem}\label{The4.1}
Suppose that the set  $ S = Sol(C, f) \cap Fix(T)$ is nonempty, the bifunction $f$ satisfies assumptions (${A}_{1bis}$), $(A_2)$, $(A_3)$ on $C$, and (${A}_{4bis}$) on $\Omega$, the mapping $T$ satisfies assumption $(A_5)$. Then the sequences $\{x^k\}$, $\{u^k\}$, $\{v^k\}$ generalized by Algorithm 3 converge weakly to $x^* \in S$, where $x^* = \lim_{k \to \infty} P_S (x^k)$.
      \end{theorem}
{\bf Proof.}
Take any $q \in S$. Since $\gamma_k \in [\b{$\gamma$}, \bar{\gamma}] \subset (0, 2)$, we deduce from Lemma~\ref{Lem4.1} that
\begin{equation}\label{4.3}
   \|u^k - q \| \leq \|x^k - q \|.
\end{equation}
By the same argument as in the proof of Theorem 1, we have
\begin{equation}\label{4.4}
\| v^k - q \|  \leq  \| x^k - q \|,
\end{equation}
and
\begin{equation}\notag
\| x^{k+1} - q \|  \leq  \| x^k - q \|.
\end{equation}
Therefore
\begin{equation}\label{4.6}
\lim_{k \to \infty} \| x^k - q \| = \tau.
\end{equation}
 Consequently, $\{x^k\}$ is bounded. Together with (\ref{4.3}), (\ref{4.4}), one has $\{u^k\}$, $\{v^k\}$ are also bounded.\\
Since
\begin{equation}\notag
\begin{aligned}
\|x^{k+1} - q \|^2 
                                      & \leq \beta_k\|v^{k} - q \|^2 + (1 - \beta_k)\|u^k - q \|^2.     \\
\end{aligned}
\end{equation}
In view of (\ref{4.4}) and Lemma~\ref{Lem4.1}, yields
\begin{equation*}
\| x^{k+1} - q \|^2 \leq   \|x^k - q\|^2 - (1 - \beta_k) \gamma_k(2-\gamma_k)(\sigma_k\|w^k\|)^2.
\end{equation*}
Therefore
\begin{equation}\label{4.11}
(1 - \beta_k) \gamma_k(2 - \gamma_k)(\sigma_k \|w^k\|)^2  \leq \big ( \|x^k - q \| - \| x^{k+1} - q \| \big ) \big( \|x^k - q \| + \| u^k - q \| \big ).
\end{equation}
Because $0 < 1 -\bar{\beta} \leq 1- \beta_k$; $\gamma_k \in [\b{$\gamma$}, \bar{\gamma}] \subset (0, 2)$, and (\ref{4.6}), we obtain from (\ref{4.11}) that
\begin{equation}\label{4.12}
\lim_{k \to \infty} \sigma_k\|w^k \| = 0.
\end{equation}
Since $u^k = P_C(x^k - \gamma_k\sigma_k w^k)$, we have
$$ \|u^k - x^k\| \leq  \gamma_k\sigma_k \|w^k\|.$$
Combining with (\ref{4.12}) we get
\begin{equation}\label{4.14}
\lim_{k \to \infty} \|u^k - x^k \| = 0.
\end{equation}
Arguing  similarly as in the proof of Theorem 1, we have
\begin{equation}\label{l4.18}
\lim_{k\to \infty} \|Tu^k - u^k \| = 0.
\end{equation}
Since $\{x^k\}$ is bounded, by Lemma~\ref{Lem2.6}, $\{y^k\}$ is bounded, consequently $\{z^k\}$ is bounded. From Lemma~\ref{Lem2.2}, $\{w^k\}$ is bounded.
In view of  (\ref{4.12}) yields
\begin{equation}\label{4.16}
  \lim_{k \to \infty} f(z^k, x^k) = \lim_{k \to \infty} [\sigma_k\|w^k\|] \|w^k\| = 0.
\end{equation}
Moreover
\begin{equation*}
\begin{aligned}
0 = f(z^k, z^k) & = f(z^k, (1 - \eta_{k}) x^k + \eta_k y^k) \\
                        & \leq (1 - \eta_k)f(z^k, x^k) + \eta_k f(z^k, y^k),
\end{aligned}
\end{equation*}
so, we get from (\ref{4.1}) that
\begin{equation*}
\begin{aligned}
f(z^k, x^k) & \geq \eta_k [f(z^k, x^k) - f(z^k, y^k)]\\
                  &\geq  \frac{\mu}{2\rho_k} \eta_k \|x^k - y^k \|^2.
\end{aligned}
\end{equation*}
In view of (\ref{4.16}) one has
\begin{equation}\label{4.17}
 \lim_{k \to \infty} \eta_k\| x^k - y^k\|^2 = 0.
\end{equation}
Suppose that $\{x^{k_i}\} \subset \{x^k\}$ and $x^{k_i} \rightharpoonup x^*$ as $i \to \infty$. From (\ref{4.17}) we get
\begin{equation}\label{4.18}
 \lim_{i \to \infty} \eta_{k_i}\| x^{k_i} - y^{k_i}\|^2 = 0.
\end{equation}
We now consider two distinct cases:

{\it Case 1.} $\lim\sup_{i \to \infty}\eta_{k_i} > 0$. \\
 In this case, there exist $\bar{\eta} > 0$ and a subsequence of $\eta_{{k}_i} $, denoted again by $\eta_{{k}_i}$ such that  $ \eta_{{k}_i} > \bar{\eta }, \ \forall i \geq i_0 $, for some $i_0 \geq 0$. Using this fact  and (\ref{4.18}),  one has
\begin{equation}\label{4.19}
\lim_{i \to \infty}{\|x^{{k}_i}-y^{{k}_i}\|} = 0.
\end{equation}
 Because  $ x^k \rightharpoonup x^*$, and (\ref{4.19}), it implies  that $ y^{k_i} \rightharpoonup x^*$ as $i \to \infty$. \\
From assertation (i) of Lemma~\ref{Lem3.1} we get
\begin{equation*}
\rho_{k_i}\big[ f(x^{k_i}, y) - f(x^{k_i}, y^{k_i}) \big] \geq  \langle x^{k_i} - y^{k_i} , y - y^{k_i} \rangle, \ \forall y \in C.
\end{equation*}
So
\begin{equation}\label{4.20}
\rho_{k_i}\big[ f(x^{k_i}, y) - f(x^{k_i}, y^{k_i}) \big] \geq   -\|x^{k_i} - y^{k_i}\| \| y - y^{k_i}\|.
\end{equation}
Letting $i \to \infty$, by the continuity of $f$ and  (\ref{4.19}), we obtain in the limit from (\ref{4.20}) that
$$ f(x^*, y) - f(x^*, x^*) \geq 0. $$
Hence
$$ f(x^*,y)  \geq 0, \ \forall y \in C, $$
which  implies that $x^*$ is a solution of EP($C, f).$

{\it Case 2.} \   $\lim_{i \to \infty}{\eta_{k_i}} = 0$.\\
From the boundedness of $\{y^{k_i}\}$, without loss of generality we may assume that $ y^{k_i} \rightharpoonup  \bar{y}$ as $i \to \infty$. \\
 Replacing $y$ by $x^{k_i}$ in ($i$) of Lemma~\ref{Lem3.1} we obtain
\begin{equation}\label{4.21}
f(x^{k_i}, y^{k_i}) \leq -\frac{1}{\rho_{k_i}} \| y^{k_i} - x^{k_i} \|^2.
\end{equation}
In the other hand, by the Armijo linesearch rule (\ref{4.1}), for $m_{k_i} - 1$, one has
\begin{equation}\notag
 f(z^{k_i, m_{k_i} - 1}, x^{k_i}) - f(z^{k_i, m_{k_i} - 1}, y^{k_i})  < \frac{\mu}{2\rho_{k_i}} \| y^{k_i}-x^{k_i}\|^2.
\end{equation}
Combining with (\ref{4.21}) we get
\begin{equation}\label{4.23}
f(x^{k_i}, y^{k_i}) \leq -\frac{1}{\rho_{k_i}} \| y^{k_i} - x^{k_i} \|^2 \leq \frac{2}{\mu} \big[f(z^{k_i, m_{k_i} - 1}, y^{k_i}) - f(z^{k_i, m_{k_i} - 1}, x^{k_i}) \big].
\end{equation}
 By the algorithm, we have $z^{k_i, m_{k_i} - 1} = (1-\eta^{m_{k_i} - 1})x^{k_i} + \eta^{m_{k_i} - 1}y^{k_i}$, $\eta^{k_i, m_{k_i} - 1} \to 0$  and $x^{k_i} $ converges weakly to $x^*$, $y^{k_i}$ converges weakly to $\bar{y}$, it implies that $z^{k_i, m_{k_i} - 1} \rightharpoonup x^*$ as $i \to \infty$. In addition $\{\frac{1}{\rho_{k_i}}\|y^{k_i} - x^{k_i}\|^2\}$ is bounded, without loss of generality, we may assume that $\lim_{i \to +\infty}\frac{1}{\rho_{k_i}}\|y^{k_i} - x^{k_i}\|^2$ exists. Therefore, we get in the limit from (\ref{4.23}) that
\begin{equation*}
f(x^*, \bar{y}) \leq -  \lim_{i \to +\infty}\frac{1}{\rho_{k_i}}\|y^{k_i} - x^{k_i}\|^2 \leq \frac{2}{\mu}f(x^*, \bar{y}).
\end{equation*}
So, $f(x^*, \bar{y}) = 0$ and $\lim_{i \to +\infty}\|y^{k_i} - x^{k_i}\|^2 = 0$. By the Case 1,
it is immediate that $x^*$ is a solution of EP($C, f$).\\
Moreover, from (\ref{4.14}) and (\ref{l4.18}), we have $u^{k_i} \rightharpoonup x^*$ and $\lim_{i \to \infty}\|Tu^{k_i} - u^{k_i}\| = 0$. By Lemma~\ref{L2.7}, $I - T$ is demiclosed at zero, hence $Tx^* = x^*$, i.e., $x^* \in Fix(T)$.\\
Therefore $x^* \in S$.\\
The rest of the proof can be done by the same way as before. \hfill$\Box$ \\


\section{Numerical example}
To illustrate the proposed algorithms, we consider a problem by taking
\begin{align*}
f(x,y) & =  (Px + Qy +r)^T(y - x), \\
Tx  & =  (I + U)^{-1}x,
\end{align*}
where $P = (p_{ij})_{n \times n}$, $Q = (q_{ij})_{n \times n}$, $U = (u_{ij})_{n \times n}$ are $n \times n$ symmetric positive semidefinite matrices such that $P-Q$ is also positive semidefinite and $r \in \mathbb{R}^n.$ The bifunction $f$ has the form of the one arising from a Nash-Cournot oligopolistic electricity market equilibrium model \cite{CKK} and that $f$ is convex in second variable, Lipschitz-type continuous with constants $L_1 = L_2 = \frac{1}{2}\|P-Q\|$. Because   $P-Q$ is positive semidefinite matrix, $f$ is monotone \cite{MQH}. It can be seen that the set of fixed points of mapping $T$ is the solution set of the equation $Ux = 0.$ In order to ensure that the intersection of the fixed points of the mapping $T$ and the solution set of EP($C, f$) is nonempty, we futher assume that the constraint set $C$ contains the original, $r = 0,$ and $U$ is a diagonal matrix such that  $u_{ii} > 0$, forall $i \in I_0$ and $u_{ii} = 0$, forall $i \not\in I_0$, for some index set $I_0 \subset \{1, 2, ..., n-1, n \}.$

We tested proposed algorithms for this example in which $C$ is the box $C = \prod_{i=1}^{n}[-10, 10],$ $P, \ Q, \ U$ are matrices of the form $A^TA$ with $A = (a_{ij})_{n \times n}$ being randomly generated in the interval $[-5, 5]$, starting point $x^0$ is randomly generated in $[-10, 10]$ and the parameters: $\alpha_0 = \beta_0 = \frac{1}{2}, \alpha_k = 1 - \frac{1}{k+2}, \beta_k = \frac{1}{2} + \frac{1}{k+3} $, and   $\rho_k  = \frac{0.5}{\|P - Q\|}$ in Algorithm 2; $\eta = 0.98$, $\mu = 0.4$, $\rho_k = 0.5, \gamma_k = 1$ in Algorithm 3.

We implement Algorithm 2 and Algorithm 3 for this problem in Matlab R2013 running on a Desktop with Intel(R) Core(TM) 2Duo CPU E8400 3GHz, and 3GB Ram.  To terminate the Algorithms, we use the stopping criteria $\|x^{k+1} - x^k\| < \epsilon$ with a tolerance $\epsilon = 10^{-6}$.

To compare with algorithms proposed in \cite{DK}, we also report the results computed with Algorithm 1 and Algorithm 2 in \cite{DK} for this problem with this data and a tolerance  $\epsilon = 10^{-3}$.

The computation results on Algorithm 2 and Algorithm 3 are reported in Table~\ref{Tab1} and Table~\ref{Tab2}, and the results on Algorithm 1 and Algorithm 2 in \cite{DK} are reported in Table~3 and Table~4, respectively, where
\begin{itemize}
\item[] {\it N.P}: the number of the tested problems; \\
\item[] {\it Average Times}: the average CPU-computation times (in second);\\
\item[] {\it Average Iteration}: the average number of iterations.
\end{itemize}
\begin{table}[!ht]
\caption{Results computed with Algorithm 2}\label{Tab1}
\renewcommand{\arraystretch}{1.6}
\begin{tabular}{|c|c|c|c|}
\hline
\ \ \ \  N.P  \ \ \ \   & \ \ \ \  Size (n)   \ \ \  \ & \ \ \ \ \ Average Times  \ \  \ \ &  \ \ \  \ Average Iterations \ \ \ \   \\
\hline
  10& 5  \ & 2.3766           & \ \ \ \  152    \ \ \ \        \\
\hline
10 & 10  &4.2141 & 223 \\
\hline
  10 & 20 & 6.7813         & 457         \\
\hline
  10 &  30 &  10.5266 & 515 \\
\hline
10 & 50    &           17.4891& 567    \\
\hline
10 & 100    &           29.2406& 674    \\
\hline
\end{tabular}
  \end{table}

\begin{table}[!ht]
\caption{Results computed with Algorithm 3}\label{Tab2}
\renewcommand{\arraystretch}{1.6}
\begin{tabular}{|c|c|c|c|}
\hline
\ \ \ \ N.P   \ \ \ \  &   \  \ \ \  Size (n)    \ \ \ \    &  \ \ \ \  Average Times  \ \ \ \  &  \ \ \ \  Average Iterations  \ \ \ \   \\
\hline
  10& 5  \ &2.4656          & \ \ \ \  99    \ \ \ \         \\
\hline
10 & 10  &4.1422 & 132 \\
\hline
  10 & 20 & 6.6375          & 164          \\
\hline
  10 &  30 &  8.0672 & 170 \\
\hline
10 & 50    &            11.8828& 192   \\
\hline
10 & 100    &            21.4953& 210   \\
\hline
\end{tabular}
  \end{table}

\begin{table}[!ht]
\caption{Results computed with Algorithm 1 in \cite{DK}}\label{Tab3}
\renewcommand{\arraystretch}{1.6}
\begin{tabular}{|c|c|c|c|}
\hline
 \ \ \ \   N.P \ \ \ \  &  \  \ \ \  Size (n)   \  \ \ \   & \ \ \ \ Average Times \ \ \ \ & \ \ \ \ Average Iterations  \ \ \ \   \\
\hline
  10& 5  \ &23.2484           & \ \ \ \  826    \ \ \ \        \\
\hline
10 & 10  &34.7438 & 1445 \\
\hline
  10 & 20 & 87.1016         & 2346         \\
\hline
  10 &  30 &  157.5781 & 2715 \\
\hline
10 & 50    &            255.4578& 3839    \\       
\hline
\end{tabular}
  \end{table}

\begin{table}
\caption{Results computed with Algorithm 2 in \cite{DK}}\label{Tab4}
\renewcommand{\arraystretch}{1.6}
\begin{tabular}{|c|c|c|c|}
\hline
 \ \ \ \   N.P \ \ \ \  & \  \ \ \  Size (n) \  \ \ \   &  \ \ \ \ Average Times \ \ \ \ & \ \ \ \ Average Iterations \ \ \ \   \\
\hline
  10& 5  \ &38.5938           & \ \ \ \ 904    \ \ \ \        \\
\hline
10 & 10  &106.3172 & 2242 \\
\hline
  10 & 20 & 163.1266          & 3050           \\
\hline
  10 &  30 &  250.9313 & 3001 \\
\hline
10 & 50    &            359.1094 & 3592   \\       
\hline
\end{tabular}
  \end{table}

From the computed results reported in these tables, we can see that the computation times and the number of iterations computed by weak convergence algorithms are much less than that computed by strong convergence algorithms, especially when the dimension of space is large. 

\section{Conclusions}
\small
We have introduced three iterative methods for finding a common point of the set of fixed points of  a symmetric generalized hybrid mapping  and the solution set of  an equilibrium problem  in a real Hilbert space. The basic iterations used in this paper are Ishikawa's process combining with the proximal point algorithm or the extragradient algorithm with or without the incorporation of a linesearch procedure. Then we have proved that the iterative sequences generated by each method converge weakly to a solution of this problem, a numerical example is also provided.\\

\noindent{\bf Acknowledgements} \\
The authors would like to thank the referees very much for their constructive comments and suggestions, especially on the presenting and the structure of the early version of their paper which helped them very much in revising the paper. 


\end{document}